\documentclass{article}
\usepackage{amsmath,amsfonts,amssymb,amsthm}
\usepackage{cases}
\usepackage{color}
\usepackage{xcolor}
\usepackage{colortbl}
\usepackage{graphicx}
\usepackage{hyperref}

\usepackage{cite}

\newcommand{\hl}[1]{{\color{black} #1}}

\usepackage{todonotes}

\allowdisplaybreaks
\numberwithin{equation}{section}

\newtheorem{theorem}{Theorem}

\numberwithin{theorem}{section}
\numberwithin{lemma}{section}
\numberwithin{corollary}{section}
\numberwithin{proposition}{section}
\numberwithin{definition}{section}
\numberwithin{remark}{section}
\numberwithin{example}{section}
\numberwithin{equation}{section}

\date{}

\begin{document}

\title{\bf On the cross-correlation of Golomb Costas permutations}
\author{\begin{tabular}[t]{c@{\extracolsep{0em}}c}
{\large Huaning Liu$^{1}$ and Arne Winterhof$^{2}$  }  \\ \\
  {\normalsize {$^{1}$Research Center for Number Theory and Its Applications }} \\
  {\normalsize {School of Mathematics, Northwest University }} \\
   {\normalsize {Xi'an 710127, China }} \\
  {\normalsize {E-mail: hnliu@nwu.edu.cn }} \\
  {\normalsize {$^{2}$Johann Radon Institute for Computational and Applied Mathematics }} \\
  {\normalsize {Austrian Academy of Sciences }} \\
   {\normalsize {Altenbergerstr.\ 69, 4040 Linz, Austria }} \\
  {\normalsize {E-mail: arne.winterhof@oeaw.ac.at }}
  \end{tabular}}

\maketitle

\begin{abstract}
In the most interesting case of {\em safe prime powers} $q$, G\'omez and Winterhof showed that a subfamily of the family of Golomb Costas permutations of $\{1,2,\ldots,q-2\}$ of size $\varphi(q-1)$ has maximal cross-correlation
of order of magnitude at most $q^{1/2}$. In this paper we study a larger family of Golomb Costas permutations and prove a weaker bound on its maximal cross-correlation.
Considering the whole family of Golomb Costas permutations we show that large cross-correlations are very rare. Finally, we collect several conditions for a small cross-correlation of two Costas permutations.
Our main tools are the Weil bound and the Szemer\'edi-Trotter theorem for finite fields.
\\

\\

\end{abstract}

\section{Introduction}

For a positive integer $n$, let $f$ be a permutation of $\{1,2,\ldots,n\}$ satisfying
$$
f(i+k)-f(i)\neq f(j+k)-f(j)
$$
for any integers $1\leq k\leq n-2$ and $1\leq i<j\leq n-k$. Such a permutation is called
a {\em Costas permutation } of $\{1,2,\ldots,n\}$. The binary $n\times n$-matrix $A=(a_{ij})_{i,j=1}^n$ defined by $a_{ij}=1$ if and only if $f(i)=j$ is called a {\em Costas array}. Costas arrays have many important applications in
radar and sonar, see for example \cite{GolombT1982,HJ2020,Jungnickel1993}.

The {\em cross-correlation } $C_{f_1,f_2}(u,v)$ between two mappings
$$
f_1, f_2: \ \{1,2,\ldots,n\}\longrightarrow \{1,2,\ldots,n\}
$$
at $(u,v)\in\mathbb{Z}^2$, $1-n\leq u,v\leq n-1$, is the number of solutions
$$
x\in\left\{\max\{1,1-u\}, \ldots, \min\{n,n-u\}\right\}
$$
of the equation $f_1(x)+v=f_2(x+u)$. For a family $\mathcal{F}$ of Costas permutations of $\{1,2,\ldots,n\}$ with $n\ge 2$, the {\em
maximal cross-correlation} $C(\mathcal{F})$ is
$$
C(\mathcal{F})=\max_{1-n\le u,v \le n-1}\mathop{\max_{f_1,f_2\in\mathcal{F}}}_{f_1\neq f_2}
C_{f_1,f_2}(u,v).
$$
In multi-user systems it is desirable that the signals have low cross-correlation, in order to minimize cross-talk. While Costas arrays guarantee ideal autocorrelation by definition, they typically exhibit poor cross-correlation. Families of Costas arrays with low cross-correlation are, therefore, particularly suitable for such applications, see for example \cite{DGHR08}.

For a prime power $q>2$ and primitive elements $g_1$ and $g_2$ of the finite field $\mathbb{F}_q$,
let $\pi_{g_1,g_2}$ be the permutation of $\{1,\ldots,q-2\}$ defined by
\begin{equation}\label{pi}
\pi_{g_1,g_2}(x)=y \qquad \hbox{if and only if} \qquad g_1^x+g_2^y=1.
\end{equation}
This construction of Costas permutations is due to Golomb, see
\cite{CohenM1991,DrakakisGRST2011,GolombT1984,Jungnickel1993,MorenoS1990}.

Let $\mathcal{U}_q$ denote the set of primitive elements of $\mathbb{F}_q$ and define
\begin{eqnarray*}
\mathcal{G}_q&=&\left\{\pi_{g_1, g_2}: \ g_1\in \mathcal{U}_q\right\}  \quad \mbox{for fixed }g_2\in \mathcal{U}_q\\
\mbox{and}\quad \mathcal{L}_q&=&\left\{\pi_{g_1, g_2}: \ g_1\in \mathcal{U}_q, \ g_2\in \mathcal{U}_q\right\}.
\end{eqnarray*}
It is easy to see that for $q=p^w$ with $p$ a prime and $w$ a positive integer we have
\begin{equation}\label{equiv}
\pi_{g_1,g_2}=\pi_{g_3,g_4}\quad \mbox{if and only if}\quad (g_3,g_4)=(g_1^{p^j},g_2^{p^j}) \quad \mbox{for some }j\in \{0,1,\ldots,w-1\}.
\end{equation}
For the convenience of the reader we give a short proof of $(\ref{equiv})$ in Section~\ref{equivsec}.
We call $(g_1,g_2)$ and $(g_3,g_4)$ {\em conjugates} if $(g_3,g_4)=(g_1^{p^j},g_2^{p^j})$ for some $j\in \{0,1,\ldots,w-1\}$.
So we have
$$
|\mathcal{G}_q|=\varphi(q-1)\gg \frac{q}{\log\log q}\qquad \mbox{and}\qquad |\mathcal{L}_q|=\frac{\varphi(q-1)^2}{w}\gg \frac{q^2}{w(\log\log q)^2},
$$
where $\varphi$ denotes Euler's totient function and $F(N)\gg G(N)$ if and only if $|G(N)|\le c F(N)$ for some absolute constant $c>0$.


For a prime power $q\geq 4$, let $t$ be the smallest prime divisor
of $\frac{q-1}{2}$ if $q$ is odd and of $q-1$ if $q$ is even. Then the maximal
cross-correlation $C(\mathcal{G}_q)$ of the family of
Golomb Costas permutations $\mathcal{G}_q$ of $\{1,\ldots,q-2\}$
satisfies
\begin{equation}\label{CGq}
C(\mathcal{G}_q)\left\{\begin{array}{ll}
\le 1+\left\lfloor q^{1/2}\right\rfloor & \hbox{if $q$ is odd and $t=\frac{q-1}{2}$ or
$q$ is even and $t=q-1$}, \\
=\frac{q-1}{t}-1 & \hbox{otherwise}.
\end{array}\right.
\end{equation}
The tables of \cite{DrakakisGRST2011} show that $C(\mathcal{L}_q)$ is larger than $C(\mathcal{G}_q)$ for some small values of
$q$.  However, for
all prime values of $q$ with $61 \le q \le 271$ and all strict prime powers $25 \le q \le 343$,
the bound $(\ref{CGq})$ is also valid for $C(\mathcal{L}_q)$. It remains an open problem to
prove the conjecture that this bound holds for $C(\mathcal{L}_q)$ up to a few exceptions of small $q\le 59$.

Note that $C(\mathcal{G}_q)$ does not depend on the choice of $g_2$ in the definition of $\mathcal{G}_q$:
Note that
$$C(\mathcal{L}_q)=\max_{0\le u,v\le q-3} \max_{f_1,f_2\in \mathcal{L}_q\atop f_1\not=f_2}C_{f_1,f_2}(u,v)$$
since
$$C_{f_1,f_2}(-u,v)=C_{f_2,f_1}(u,-v)\quad \mbox{for }1\le u\le q-3$$ and
$$C_{\pi_{g_1,g_2},\pi_{g_3,g_4}}(u,-v)=C_{\pi_{g_1,g_2^{-1}},\pi_{g_3,g_4^{-1}}}(u,v)\quad \mbox{for }1\le v\le q-3.$$
Hence, we may restrict ourselves to pairs $(u,v)$ with $0\le u,v\le q-3$.
Now take another primitive element $g_4=g_2^j$ with $\gcd(j,q-1)=1$ and put $v'=j^{-1}v\bmod (q-1)$.
Then we have
$$\pi_{g_1,g_4}(x)=j\pi_{g_1,g_2}(x)\bmod (q-1)$$
and
$$\pi_{g_1,g_4}(x)+v=\pi_{g_3,g_4}(x+u)$$
is equivalent to
$$\pi_{g_1,g_2}(x)+v'=\pi_{g_3,g_2}(x+u).$$

In particular, if $t=\frac{q-1}{2}$ or $t=q-1$, respectively, we get a bound of order of magnitude $q^{1/2}$ as well as
$$|\mathcal{L}_q|=\left\{\begin{array}{ll}\frac{(q-3)^2}{4w} & \hbox{if $q$ is odd},\\
\frac{(q-2)^2}{w} & \hbox{if $q$ is even}.\end{array}\right.$$
In this case we have either
\begin{itemize}
    \item $q$ is a {\em safe prime}, that is, $\frac{q-1}{2}$ is a ({\em Sophie-Germain}) prime,
    \item $q-1$ is a {\em Mersenne prime}, that is, $q=2^w$ for some prime $w$ such that $q-1$ is a prime,
    \item or $q$ is a {\em strict safe prime power}, that is, $q=3^w$ for some prime $w>2$ such that $\frac{q-1}{2}$ is a prime by \cite[Lemma~1]{DrakakisGRST2011}.
\end{itemize}
We focus on these three cases \hl{since otherwise the bound on $C({\cal G}_q)$ is larger as well as the size of ${\cal L}_q$ is smaller (in terms of $q$)}. By a slight abuse of notation we call all such $q$ {\em safe prime powers}. The first such $q$ are
$$q=4,5,7,8,11,23,27,32,47,59,83,107,128,167,179,227,263.$$

The maximal cross-correlation
of $\mathcal{L}_q$ can be bounded by the number of solutions of certain polynomial equations, see \cite{DrakakisGRST2011},
$$
C(\mathcal{L}_q)\leq \mathop{\max_{\alpha,\beta\in\mathbb{F}_q^{*}}}_{1\leq r,s\leq q-2, \ \gcd(rs, q-1)=1}\left|\left\{x\in\mathbb{F}_q: \ \alpha(1-x)^{r}+\beta x^{s}-1{\color{blue}=0}\right\}\right|,
$$
However,  without restrictions on $s$ and $r$ a bound on the number of solutions of
the equation
$$
\alpha(1-x)^{r}+\beta x^{s}-1=0
$$
seems to be out of reach with currently known methods.
Still finding larger subfamilies of $\mathcal{L}_q$ than $\mathcal{G}_q$ with non-trivial cross-correlation of size at most $q^{1-\varepsilon}$ for some $\varepsilon>0$, is an important open problem. It is expected that the cross-correlation increases with the size of the subfamily.

Our first result is a weaker bound than $(\ref{CGq})$ but for a larger subfamily of the family of Golomb Costas permutations.

{\color{blue}\begin{theorem}\label{Theorem:Golomb-5}
Let $q>7$ be a safe prime power. Fix $0 <\delta <1/2$.
Then there exists a subfamily $\mathcal{L}_{q,\delta}$ of $\mathcal{L}_q$ of size $|\mathcal{L}_{q,\delta}| \gg \delta q\log q$ with cross-correlation $C(\mathcal{L}_{q,\delta})=O\left(q^{1/2+\delta}\right)$.
\end{theorem}
}

We prove \hl{(a slightly more precise version of)} Theorem~\ref{Theorem:Golomb-5} in Section~\ref{proof11} \hl{in a constructive way}.

Write
$$g_3=g_1^r\quad \mbox{and}\quad g_4=g_2^s \quad \mbox{with}\quad  1\leq r,s<q-1 \mbox{ and }\gcd(rs, q-1)=1.$$
Theorem~\ref{Theorem:Golomb-5} is based on
\begin{equation}\label{CLqdelta}
C_{\pi_{g_1,g_2},\pi_{g_3,g_4}}(u,v)\le 1+\min\{s,s^{-1}\}q^{1/2}\quad \mbox{if $(g_1,g_2)$ and $(g_3,g_4)$ are not conjugates},
\end{equation}
where $s^{-1}$ denotes the inverse of $s$ modulo $q-1$, that is, $1\le s^{-1}\le q-2$,
and as the bound of \cite{Gomez-PerezW2020} on the Weil bound, which is optimal for some polynomials.

Our second result states that for any fixed pair of Golomb-Costas permutations the proportion of $(u,v)$ with large cross-correlation at $(u,v)$ is very small.

For a fixed pair $(u,v)$ we also prove a similar result on the number of quadruples $(g_1,g_2,g_3,g_4)$ with large cross-correlation $C_{\pi_{g_1,g_2},\pi_{g_3,g_4}}(u,v)$.

More precisely, the following theorem suggests that
the number of possibly large cross-correlations is small in both cases, for fixed $(g_1,g_2,g_3,g_3)$ and for fixed $(u,v)$, respectively. The proof depends on an analog of the Szemer\'edi-Trotter theorem for finite fields, that is, a point-line incidence bound. It seems that this technique has not been used for studying cross-correlations of Costas arrays yet.

\begin{theorem}\label{Theorem:Golomb-3} 
1. Let $\pi_{g_1,g_2}, \pi_{g_3,g_4}\in\mathcal{L}_q$,
$B>0$, ${\cal S}\subset \{0,1,\ldots,q-3\}^2$,
and let
$N_{q,B, {\cal S}}$ be the number of $(u,v)\in {\cal S}$ with
$$
C_{\pi_{g_1,g_2},\pi_{g_3,g_4}}(u,v)\ge B.
$$
Then we have for any $q$,
$$
N_{q,B,{\cal S}}\ll \frac{1}{B} \left\{ \begin{array}{cc} |{\cal S}|^{1/2}q, & |{\cal S}|\ge q,\\
|{\cal S}|q^{1/2}, & q^{1/2}\le |{\cal S}|<q,\\
q, & |{\cal S}| <q^{1/2}.
\end{array}\right.
$$
and for prime $q$ the improvement,
$$N_{q,B,{\cal S}}\ll \frac{(|{\cal S}|q)^{11/15}}{B}\quad \mbox{if }  q^{7/8} < |{\cal S}|< q^{8/7},\quad q \mbox{ prime}.$$
2. Let $u, v\in\{0,1,\ldots, q-3\}$, $B>0$ and let $M_{q,B}$ be the number of
$(\pi_{g_1,g_2}, \pi_{g_3,g_4})\in\mathcal{L}_q^{2}$ such that
$$
C_{\pi_{g_1,g_2},\pi_{g_3,g_4}}(u,v)\ge B.
$$
Then we have
$$
M_{q,B}<  \frac{\tau(q-1)\varphi(q-1)^3(q-1)}{B}\ll \frac{q^{4+o(1)}}{B},
$$
where $\tau(n)$ is the number of positive divisors of $n$ and $f(n)=o(1)$ if $\lim\limits_{n\rightarrow \infty} f(n)=0$.
\end{theorem}
Although for multi-user systems it would be nice to have small cross-correlation for all $(u,v)$ it may be still sufficient to allow some exceptions which appear with rather small probability, which is made explicit in the first part of Theorem~\ref{Theorem:Golomb-3}. Similarly, cross-talk is limited even if, for fixed $(u,v)$, exceptions of pairs of Costas arrays with large cross-correlation occur only with small probability, which motivates the second part.


We prove Theorem~\ref{Theorem:Golomb-3} in Section~\ref{secadd}.

Finally, in Section~\ref{Sec:improve} we show that in (\ref{CLqdelta}) we can substitute $\min\{s,s^{-1}\}$ by the minimum of
\begin{equation}\label{minset} s, s^{-1}, r, r^{-1}, q-s, q-r, q-r^{-1},q-s^{-1}.
\end{equation}
For odd $q$ we can also estimate the cross-correlation by
$$1+2\min\left\{\left|k-\frac{q-1}{2}\right|:k\in \{s,r,s^{-1},r^{-1},q-s,q-r,q-s^{-1},q-r^{-1}\right\}q^{1/2}.$$

Alternatively, we can estimate the cross-correlation by $\max\{r,s\}$ and substitute either $r$ by $q-1-r$, $s$ by $q-1-s$ or both.
For odd $q$ there are also obvious improvements for certain pairs of $(r,s)$ if $r$, $s$ or both is close to $\frac{q-1}{2}$.

\section{Preliminary results}

\subsection{Proof of $(\ref{equiv})$}
\label{equivsec}

With $g_3=g_1^r$ and $g_4=g_2^s$ for some $r$ and $s$ with $\gcd(rs,q-1)=1$ and $1\le r,s\le q-2$, the permutations $\pi_{g_1,g_2}$ and $\pi_{g_3,g_4}$ coincide if and only if
$$(1-g_1^x)^s=g_2^{s\pi_{g_1,g_2}(x)}=g_4^{\pi_{g_3,g_4}(x)}=1-g_1^{rx}\quad \mbox{for }x=1,2,\ldots,q-2,$$
by $(\ref{pi})$,
that is, the polynomial equation
\begin{equation}\label{poly1}
(1-X)^s=1-X^r
\end{equation}
of degree at most $q-2$ has at least $q=p^w$ different solutions, namely, $0,1,g_1,g_1^2,\ldots,g_1^{q-2}$ and thus the polynomials on both sides of $(\ref{poly1})$ coincide.
In particular, we get $r=s$. Write $s=s_0p^{j}$ with $1\le s_0<p$ and $0\le j<w$. Substitute
$X^{p^{j}}$ by $Y$ to get
\begin{equation}\label{poly2}(
1-Y)^{s_0}=(1-X)^s=(1-X^s)=1-Y^{s_0}.
\end{equation}
Since $\xi\mapsto \xi^{p^{j}}$ is a bijection on $\mathbb{F}_q$ this equation of degree $s_0<p$ has still $q\ge p$ solutions and the polynomials on both sides of $(\ref{poly2})$ coincide.
The multiplicity of $1$ as a zero of the left hand side is $s_0$ and of the right hand side $1$ since the derivative $(1-Y^{s_0})'=-s_0Y^{s_0-1}$ does not vanish at $1$.
Hence, $s_0=1$ and thus $r=s=p^{j}$. \hfill $\Box$

\subsection{The Weil bound}

Let $\chi$ be a multiplicative character of order $d>1$ of the finite field
$\mathbb{F}_q$ of $q$ elements and
let $f(X)\in \mathbb{F}_q[X]$ be a polynomial which is not
a $d$-th power of a polynomial over the algebraic closure $\overline{\mathbb{F}}_q$ of $\mathbb{F}_q$. From Theorem 5.41 of \cite{LidlN1997} or Theorem 2C in \cite{Schmidt1976} we have
\begin{eqnarray}\label{Equation:character sums}
\left|\sum_{x\in \mathbb{F}_q}\chi\left(f(x)\right)\right|\leq
(m-1)q^{1/2},
\end{eqnarray}
where $m$ denotes the number of distinct zeros of $f(X)$ in $\overline{\mathbb{F}}_q$.
Note that the condition on $f(X)$ is satisfied if $f(X)$ has a simple root.

\subsection{Point-line incidence bound}

Given a finite set of points $P$ and a finite set of lines $L$ in the plane $\mathbb{F}_q^2$, $q=p^w$, we define
$$
I(P, L)=\left|\left\{(z, \ell)\in P\times L: \ z\in \ell\right\}\right|
$$
to be the number of incidences between $P$ and $L$. Several bounds on $I(P,L)$ can be found in the literature, see 
\cite{BKT2004,D2010,Mohammadi2020,PhuongTV2019,StevensZ2017,Vinh2011}. In particular, we have the following analog of the Szemer\'edi-Trotter theorem, see \cite{StevensZ2017} and references therein:
\begin{eqnarray}\label{Equation:incidence bound}
I(P, L)\ll \left\{ \begin{array}{ll} |L|, & |P|^2\le |L|,\\
|P||L|^{1/2}, & |P|\le |L|<|P|^2,\\
|P|^{1/2}|L|, & |P|^{1/2}\le |L|<|P|,\\
|P|, & |L|<|P|^{1/2},\end{array}\right.
\end{eqnarray}
with the following improvement in a nontrivial range if $|P|<p^{8/5}$,
\begin{eqnarray}\label{Equation:incidence bound2}
I(P, L)\ll (|P||L|)^{11/15}, \quad |P|^{7/8}< |L|<\min\{|P|^{8/7},|P|^{2/13}p^{15/13}\}, \quad q=p^w.
\end{eqnarray}

\section{Proof of Theorem \ref{Theorem:Golomb-5} }\label{proof11}

{\color{blue}
Let $q>7$ be a safe prime power and
$$
a=\left\{\begin{array}{cc}
3, & q\not=3^w,\\
5, & q=3^w.
\end{array}\right.
$$
Let $g$ be a fixed primitive element of $\mathbb{F}_q$.
For any $0\le \delta<\frac{1}{2}$ put
$$
\mathcal{A}_{q,a,\delta}=\left\{g^{a^{i}}: \ i=0,\pm 1,\pm 2,\ldots,\pm \left\lfloor\frac{\delta\log q}{2\log a}\right\rfloor
\right\}
$$
and
$$
\mathcal{L}_{q,\delta}=\left\{\pi_{g_1, g_2}: \ g_1\in \mathcal{U}_q, \ g_2\in \mathcal{A}_{q,a,\delta}\right\}.
$$ }

Note that $a$ is the smallest integer $>1$ with $\gcd(a,q-1)=1$ and $a\not=p$, where $q=p^w$.
By this choice of $a$ two different pairs $(g_1,g_2)$ and $(g_3,g_4)$ in $\mathcal{U}_q\times \mathcal{A}_{q,a,\delta}$ are not conjugates and thus $\pi_{g_1,g_2}, \pi_{g_3,g_4}\in\mathcal{L}_{q,\delta}$ are different by $(\ref{equiv})$.
{\color{blue}Hence, we get the formula for the size of $\mathcal{L}_{q,a,\delta}$ as follows
$$
|\mathcal{L}_{q,\delta}|=\varphi(q-1)\left(2\left\lfloor\frac{\delta\log q}{2\log a}\right\rfloor+1\right)=\left\{\begin{array}{ll}\frac{q-3}{2}\left(2\left\lfloor\frac{\delta\log q}{2\log a}\right\rfloor+1\right) & \hbox{if $q$ is odd},\\ (q-2)\left(2\left\lfloor\frac{\delta\log q}{2\log a}\right\rfloor+1\right) & \hbox{if $q$ is even}. \end{array}\right.
$$
}

Write again
$$g_3=g_1^r\quad \mbox{and}\quad g_4=g_2^s \quad \mbox{with}\quad  1\leq r,s<q-1 \mbox{ and }\gcd(rs, q-1)=1.$$
First we prove $(\ref{CLqdelta})$.

For $s=1$ we get $(\ref{CLqdelta})$ from $(\ref{CGq})$ and hence we may assume that $s>1$.

By \cite[Theorem~2]{DrakakisGRST2011} we have
\begin{eqnarray}\label{Ceq}
C_{\pi_{g_1,g_2},\pi_{g_3,g_4}}(u,v)&\leq&\left|\{y\in\mathbb{F}_q\setminus\{1\} : g_1^{ru}(1-y)^{r}=1-g_4^{v}y^{s}\}\right|\\
&=&\left|\{z\in\mathbb{F}_q\setminus\{g_4^v\} : g_1^{-u}(1-z)^{r^{-1}}
=1-g_4^{-vs^{-1}}z^{s^{-1}}\}\right|,\nonumber
\end{eqnarray}
where we substituted $z=g_4^{v}y^s$. 
Hence, we may restrict ourselves to the case $s\le s^{-1}$.

As in the proof of \cite[Theorem~2]{Gomez-PerezW2020} we get
\begin{eqnarray*}
&&\left|\{y\in\mathbb{F}_q\setminus\{1\} : g_1^{ru}(1-y)^{r}=1-g_4^{v}y^{s}\}\right|\\
&=&\frac{1}{q-1}\sum_{\chi}\sum_{y\in\mathbb{F}_q}
\chi(g_1^{ru}(1-y)^{r}(1-g_4^{v}y^{s})^{q-2}),
\end{eqnarray*}
where $\chi$ runs through all multiplicative characters of $\mathbb{F}_q$ and we use the convention $\chi(0)=0$.
The contribution of the principal character is at most $1$ and for the non-principal characters
from the Weil bound $(\ref{Equation:character sums})$ we get
$$
\left|\sum_{y\in\mathbb{F}_q}
\chi\left(g_1^{ru}(1-y)^{r}(1-g_4^{v}y^{s})^{q-2}\right)\right|\leq sq^{1/2},\quad s>1,
$$
and $(\ref{CLqdelta})$ follows. Note that the Weil bound is applicable since the polynomial $1-g_4^vY^s$ has a simple root~$\not=1$ since $s>1$
and we may assume $\gcd(s,p)=1$ by $(\ref{equiv})$.


Put $\theta=\left\lfloor\frac{\delta\log q}{2\log a}\right\rfloor$, $g_2=g^{a^i}$ and $g_4=g^{a^j}=g_2^s$ for some $i$ and $j$ with $0\le |i|,|j|\le \theta$. Then we have
\begin{eqnarray*}
s\equiv a^{j-i} \bmod (q-1)
\quad\mbox{ and }
s^{-1}\equiv 
a^{i-j}\bmod (q-1). 
\end{eqnarray*}
Therefore
$$C_{\pi_{g_1,g_2},\pi_{g_3,g_4}}(u,v)\leq 1+a^{|i-j|}q^{1/2}\le 1+a^{2\theta}q^{1/2}\leq 1+q^{1/2+\delta}$$
and the result follows. \hfill $\Box$


\section{Proof of Theorem \ref{Theorem:Golomb-3} }\label{secadd}

Let $\pi_{g_1,g_2}, \pi_{g_3,g_4}\in\mathcal{L}_q$ with
$g_3=g_1^r$ and $g_4=g_2^s$
with $1\le r,s<q-1$ and $\gcd(rs, q-1)=1$. Define a set~$P$ of points and a set $L$ of lines,
\begin{eqnarray*}
P&=&\left\{\left((1-x)^s, \ x^r\right): \ x\in\mathbb{F}_q\right\},\\
L&=&\left\{g_4^vy+g_1^{ru}z=1: \ (u, v)\in \cal S\right\},
\end{eqnarray*}
of sizes  $|P|=q$ and $|L|= |{\cal S}|$.
By \cite[Theorem~2]{DrakakisGRST2011} we get
\begin{eqnarray*}
BN_{q,B,{\cal S}} &\le&
\sum_{(u,v)\in {\cal S}} C_{\pi_{g_1,g_2},\pi_{g_3,g_4}}(u,v)\\
&\leq&  \sum_{(u,v)\in {\cal S}}
\left|\{x\in\mathbb{F}_q :  g_4^{v}(1-x)^s+g_1^{ru}x^r=1\}\right|  \\
&=&I(P, L)
\end{eqnarray*}
and the first result follows from $(\ref{Equation:incidence bound})$
and $(\ref{Equation:incidence bound2})$.

Moreover, for fixed $(u, v)\in\{0,1,\ldots,q-3\}^2$ we get from \cite[Theorem~2]{DrakakisGRST2011} that
\begin{eqnarray*}
B M_{q,B}&\le&\mathop{\sum_{\pi_{g_1,g_2},\pi_{g_3,g_4}\in\mathcal{L}_q}}
C_{\pi_{g_1,g_2},\pi_{g_3,g_4}}(u,v)\\
&\leq& \sum_{g_1, g_2\in \mathcal{U}_q}\mathop{\sum_{s=1}^{q-1}}_{\gcd(s, q-1)=1}\mathop{\sum_{r=1}^{q-1}}_{\gcd(r, q-1)=1}
\left|\{x\in\mathbb{F}_q^*\setminus\{1\} :g_1^{ru}x^r =1- g_2^{sv}(1-x)^s\}\right|.
\end{eqnarray*}
Now for a fixed triple $(g_1,g_2,s)$ we estimate the inner sum over $r$, that is, we have to estimate the number of pairs $(r,x)$
with
$$(g_1^ux)^r=1-(g_2^v(1-x))^s.$$
Note that $x\mapsto g_1^u x$ is a bijection of $\mathbb{F}_q^*$ and that for each divisor $d$ of $q-1$ there are $\varphi\left(\frac{q-1}{d}\right)$ different elements $x$ in $\mathbb{F}_q^*$ of the form $x=g^j$ with $j=0,1,\ldots,q-2$ and $\gcd(j,q-1)=d$, where $g$ is any fixed primitive element of $\mathbb{F}_q$.
Now assume that
$$g_1^ux=g^j\quad\mbox{with}\quad \gcd(j,q-1)=d$$
and define $k$ by
$$1-(g_2^v(1-x))^s=g^k,\quad 1\le k\le q-2.$$
We have to estimate the number of $r$ with $g^{jr}=g^k$, that is,
$$jr\equiv k\bmod q-1.$$
There is no solution if $d$ does not divide $k$ and otherwise we get
$$\frac{j}{d}r\equiv \frac{k}{d}\bmod \frac{q-1}{d}$$
which has a unique solution $r$ modulo $\frac{q-1}{d}$ and thus $d$ solutions $r$ with $1\le r\le q-1$.
Hence, the number of pairs $(r,x)$ is at most
$$\sum_{d|q-1}\varphi\left(\frac{q-1}{d}\right)d< \tau(q-1)(q-1).$$
Therefore
\begin{eqnarray*}
BM_{q,B}< \tau(q-1)(q-1)\varphi(q-1)^3
\end{eqnarray*}
and the second result follows.
\hfill $\Box$

\section{A collection of small cross-correlations}\label{Sec:improve}

Substituting $y=1-x$ in the equation
\begin{equation}\label{starteq}
    g_1^{ru}(1-y)^r=1-g_4^vy^s
\end{equation}
appearing in (\ref{Ceq}) we interchange the roles of $r$ and $s$. Hence, we may restrict ourselves to the case that the minimum of (\ref{minset}) depends on $s$ and not on $r$.



Multiplying (\ref{starteq}) by $y^{q-1-s}$ we get
$$g_1^{ru}(1-y)^ry^{q-1-s}=(y^{q-1-s}-g_4^r),\quad y\not=0.$$
The left hand side has two different zeros and the right hand side at most $q-1-s$. Dealing as in Section~\ref{proof11} and using the Weil bound we can estimate the corresponding character sums by
$(q-s)q^{1/2}$. Again as in Section~\ref{proof11} we can substitute $s$ by $s^{-1}$ in this bound.


For odd $q$ write $r=\frac{q-1}{2}+r_0$ for some $r_0$ with $|r_0|<\frac{q-1}{2}$ and recall that $x^{(q-1)/2}=1$ if $x$ is a nonzero square in $\mathbb{F}_q$ and $x^{(q-1)/2}=-1$ if $x$ is a nonzero non-square in $\mathbb{F}_q$ and we can split $(\ref{starteq})$ into two equations of the form $a(1-x)^s=1\pm bx^{r_0}$. Each of them has at most $\min\{|r_0|, q-|r_0|\}q^{1/2}$ solutions.

We can also estimate the number of solutions of $(\ref{starteq})$ with Lagrange's theorem by $\max\{r,s\}$.
Multiplying $(\ref{starteq})$ by $y^{q-1-s}$ we get a polynomial equation of degree $q-1-s+r$ with at most $q-1-s+r$ solutions.
Interchanging the roles of $r$ and $s$ and splitting the equation into two equations if either $r$ or $s$ is close to $\frac{q-1}{2}$
we get several other obvious improvements.

\section*{\bf \Large Acknowledgments}

The authors express their gratitude to the associate editor and the referees for their helpful and detailed comments.
The first author is supported by National Natural Science Foundation of China under Grant No. 12071368,
the Science and Technology Program of
Shaanxi Province of China under Grant No. 2024JC-JCQN-04, and
Shaanxi Fundamental Science Research Project for Mathematics and Physics under Grant No. 22JSY017.




\begin{thebibliography}{99}


\bibitem{BKT2004} {\rm J. Bourgain, N. Katz, T. Tao},
A sum-product estimate in finite fields, and applications.
{\em Geom. Funct. Anal.,} 14 (2004), no. 1, pp.~27--57.

\bibitem{CohenM1991} {\rm S. D. Cohen and G. L. Mullen},
Primitive elements in finite fields and Costas arrays,
{\em Applicable Algebra Eng., Commun. Comput., } 2 (1991), no. 1, pp.~45--53.

\bibitem{DGHR08}
{\rm K. Drakakis, R. Gow, J. Healy, S. Rickard},
Cross-correlation properties of Costas arrays and their images under horizontal and vertical flips.
{\em Math. Probl. Eng.} (2008), Art.\ ID 369321, 11 pp.

\bibitem{DrakakisGRST2011} {\rm K. Drakakis, R. Gow, S. Rickard, J. Sheekey and K. Taylor},
On the maximal cross-correlation of algebraically constructed Costas arrays,
{\em IEEE Trans. Inform. Theory, } 57 (2011), no. 7, pp.~4612--4621.

\bibitem{D2010} {\rm Z. Dvir},
Incidence theorems and their applications.
{\em Found. Trends Theor. Comput. Sci.,} 6 (2010), no. 4, pp.~257--393.

\bibitem{GolombT1982} {\rm S. Golomb and H. Taylor},
Two-dimensional synchronization patterns for minimum ambiguity,
{\em IEEE Trans. Inform. Theory, } 28 (1982), no. 4, pp.~600--604.

\bibitem{GolombT1984} {\rm S. Golomb and H. Taylor},
Constructions and properties of Costas arrays,
{\em Proc. IEEE, } 72 (1984), no. 9, pp.~1143--1163.

\bibitem{Gomez-PerezW2020} {\rm D. G\'{o}mez-P\'{e}rez and A. Winterhof},
A note on the cross-correlation of Costas permutations,
{\em IEEE Trans. Inform. Theory, } 66 (2020), no. 12, pp.~7724--7727.

\bibitem{HJ2020} {\rm D. Hachenberger and D. Jungnickel}, Topics in Galois fields.
Algorithms and Computation in Mathematics, 29. Springer, Cham,  2020.

\bibitem{Jungnickel1993} {\rm D. Jungnickel},
Finite Fields: Structure and Arithmetics. Mannheim, Germany: Bibliographisches Institut, 1993.

\bibitem{LidlN1997} R. Lidl and H. Niederreiter, Finite Fields (Encyclopedia of Mathematics
and Its Applications), 2nd ed. Cambridge, U.K.: Cambridge Univ. Press,
1997.

\bibitem{Mohammadi2020} {\rm A. Mohammadi},
Szemer\'{e}di-Trotter type results in arbitrary finite fields,
{\em Integers, } 20 (2020), Paper No. A7, 29 pp.

\bibitem{MorenoS1990} {\rm O. Moreno and J. Sotero},
Computational approach to conjecture a of Golomb, in
{\em Proc. 20th Southeastern Conf. Combinatorics, Graph Theory, Comput., }
Boca Raton, FL, USA, 70 (1990), pp.~7--16.

\bibitem{PhuongTV2019} {\rm N. D. Phuong, P. V. Thang and L. A. Vinh},
Incidences between planes over finite fields,
{\em Proc. Amer. Math. Soc., } 147 (2019), no. 5, pp.~2185--2196.

\bibitem{Schmidt1976} W. M. Schmidt,
Equations over finite fields. An elementary approach. Lecture Notes in Math. 536,
Springer, New York, 1976.

\bibitem{StevensZ2017} {\rm S. Stevens and F. de Zeeuw},
An improved point-line incidence bound over arbitrary fields,
{\em Bull. London Math. Soc., } 49 (2017), no. 5, pp.~842--858.

\bibitem{Vinh2011} {\rm L. A. Vinh},
The Szemer\'{e}di-Trotter type theorem and the sum-product estimate in finite fields,
{\em European J. Combin., } 32 (2011), no. 8, pp.~1177--1181.

\end{thebibliography}
\end{document}